\documentclass[11pt]{article} 
\usepackage{amssymb,latexsym}

\newtheorem{theorem}{Theorem}[section]
\newtheorem{proposition}[theorem]{Proposition}
\newtheorem{corollary}[theorem]{Corollary}

\newtheorem{lemma}[theorem]{Lemma}

\def\remark{ \textbf{Remark.} }

\def\beq{\begin{equation}}
\def\eeq{\end{equation}}

\def\A{\mathbb A}
\def\C{\mathbb C}
\def\G{\mathbb G}

\def\Z{\mathbb Z }

\def\nf{\mathfrak{n}}
\def\tf{\mathfrak{t}}

\def\P{\mathbb{P}}

\def\O{\Omega}

\def\O{{\cal O}}

\def\<{\langle\kern-.08cm\langle}
\def\>{\rangle\kern-.08cm\rangle}

\def\a{\alpha}

\def\bbf{\mathfrak{b}}

\date{}

\title{%
\textbf{Positivity in equivariant Schubert calculus}}
\author{%
William Graham${}^1$\\[.2in]
}

\begin{document}

\maketitle

\footnotetext[1]{University of Georgia, Department of Mathematics,
Boyd Graduate Studies Research Center, Athens, GA 30602}


\section{Introduction}

Let $X = G/B$ be the flag variety of a complex semisimple group $G$
with $B \supset T$ a Borel subgroup and maximal torus, respectively.
The homology $H_*(X)$ has as a basis the
fundamental classes $[X_w]$ of Schubert varieties $X_w \subset X$; if
$\{x_w\} \subset H^*(X)$ is the corresponding dual basis for
cohomology, the cup product, expressed in this basis, has nonnegative
coefficients:
\begin{equation} \label{e.1}
x_u x_v = \sum a^w_{uv} x_w 
\end{equation}
where $a^w_{uv}$ are nonnegative integers.  

The $T$-equivariant cohomology and Chow groups of the flag variety
have been described by \cite{Arabia}, \cite{Kostant-Kumar},
\cite{Brion}.  One reason to study these groups is that they provide a
way to compute the coefficients in the multiplication in ordinary
cohomology.  In addition, the equivariant groups are related to
degeneracy loci in algebraic geometry (see \cite{Fulton2},
\cite{Fulton3}, \cite{Pragacz-Ratajski}, \cite{Graham}), which in turn
are related to the double Schubert polynomials first defined in
combinatorics by \cite{Lascoux-Schutzenberger}.

Peterson \cite{Peterson} recently conjectured that the equivariant
cohomology groups of the flag variety have a positivity property
generalizing (\ref{e.1}).  The $T$-equivariant cohomology $H^*_T(X)$
is a free module over $H^*_T(pt)$ with a basis dual (in a suitable
sense; see Section \ref{s.preliminaries}) to the equivariant
fundamental classes $[X_w]_T$; again we call this basis $\{x_w\}$.
Now $H^*_T(pt)$ is isomorphic to the polynomial ring $S(\hat{T}) =
\Z[\lambda_1, \ldots, \lambda_n]$, where $\lambda_1, \ldots,
\lambda_n$ is a basis for the free abelian group $\hat{T}$ of
characters of $T$.  Let $\alpha_1, \ldots, \alpha_n$ denote the simple
roots in $\hat{T}$ (chosen so that the roots of $\bbf = \mbox{Lie }B$
are positive).  In the equivariant setting, we can again expand the
product $x_u x_v$ in the form (\ref{e.1}), but now the $a^w_{uv}$ are
in $H^*_T(pt)$ -- in other words, they are polynomials.  Peterson's
conjecture is that when each $a^w_{uv}$ is written as a sum of
monomials in the $\alpha_i$, the coefficients are all nonnegative.  In
this paper we prove the conjecture, not just for finite-dimensional
flag varieties, but in the general Ka\v{c}-Moody setting.  An
immediate corollary is a conjecture of Billey
\cite{Billey}.

The methods of this paper are those used by Kumar and Nori
\cite{Kumar-Nori}.  In that paper, the authors prove the nonnegativity
result (\ref{e.1}) in ordinary cohomology for the flag variety of a
Ka\v{c}-Moody group.  As they observe, the difficulty in proving this
result is that in the Ka\v{c}-Moody case, unlike the finite dimensional
case, the flag variety is not in general a homogeneous space.
However, it is approximated by finite-dimensional varieties, each of
which has an action of a unipotent group with finitely many orbits.
The main result of \cite{Kumar-Nori} is that for such varieties, the
cup product has nonnegative coefficients (with respect to a suitable
basis); the result for the flag variety follows.  A similar problem
arises in equivariant cohomology.  The equivariant cohomology of $X$
is by definition the cohomology of a ``mixed space'' $X_T$, which,
although infinite-dimensional, can be approximated by
finite-dimensional varieties.  As in the situation considered by Kumar
and Nori, the space $X_T$ is not a homogeneous space.  But unlike
their situation, the finite-dimensional approximations to $X_T$ do not
(as far as I know) have actions of unipotent groups with finitely many
orbits, so we cannot apply their result.  Instead, by adapting their
proof to the equivariant setting, and using a relation in equivariant
cohomology (or Chow groups) observed by Brion, we are able to deduce
an equivariant analogue of the main result of \cite{Kumar-Nori}.  The
equivariant nonnegativity result for the flag variety follows
immediately.

\medskip

{\bf Acknowledgements.}  I would like to thank Michel Brion and James
Carrell for some useful e-mail.

\section{Preliminaries} \label{s.preliminaries}

We will work with schemes over the ground field $\C$ and assume (to
freely apply the results of \cite[Ch.19]{Fulton}) that all schemes
considered admit closed embeddings into nonsingular schemes.  We use
equivariant cohomology and Borel-Moore homology with integer
coefficients as our main tools; $H_*X$ will denote the Borel-Moore
homology of $X$.  For smooth varieties, we could alternatively use
equivariant Chow groups, but for nonsmooth varieties, the Chow
``cohomology'' theory is not as well understood and for this reason we
use (equivariant) cohomology and Borel-Moore homology groups.  In this
section we recall some basic facts about these groups; for more
background, see \cite{Brion2} or \cite{Edidin-Graham1}.  We also
prove, for lack of a reference, equivariant versions of several
familiar non-equivariant results.

Let $X$ be a scheme with an action of a linear algebraic group $G$.
Let $V$ be a representation of $G$ and $U$ an open subset of $V$ such
that $G$ acts freely on $U$ and such that the (complex) codimension of
$V-U$ in $V$ is greater than $\mbox{dim }X - i/2$.  View $G$ as acting
on the right on $U$, and on the left on $X$; then $G$ acts on $U
\times X$ by $g \cdot (u,x) = (u g^{-1}, gx)$. \footnote{Alternatively,
we could let $G$ act on the left on $U$ and then take the diagonal action
on $U \times X$.}  Define $U \times^G X$ to be $(U \times X)/G$.  The
equivariant cohomology and Borel-Moore homology of $X$ are, by
definition,
\begin{eqnarray*}
H^i_G(X) &=& H^i(U \times^G X)\cr 
H^G_i(X) &=& H_{i+ 2(\dim V - \dim
G)}(U \times^G X).
\end{eqnarray*}
These groups are
independent of the choice of $V$ and $U$ provided the codimension
condition is satisfied.  For this reason we often denote $U
\times^G X$ by $X_G$ (omitting $U$ from the notation).  The quotient
$U/G$ is a finite-dimensional approximation to the classifying
space $BG$ introduced in Chow theory by Totaro \cite{Totaro}.  We
will frequently write $BG$ when we mean such a finite-dimensional
approximation.

The equivariant cohomology of a point we denote by $H^*_G$.  Both
$H^*_G(X)$ and $H^G_*(X)$ are modules for $H^*_G$.  $H^*_G(X)$ has a
natural ring structure, and $H^G_*(X)$ is a module for this ring.  Any
$G$-stable closed subvariety $Y \subset X$ has a fundamental class
$[Y]_G$ in $H^G_{2 \dim Y}(X)$ .  There is a natural map $\cap [X]_G:
H^*_G(X) \rightarrow H^G_*(X)$; if $X$ is smooth this is an
isomorphism.  In particular, we will always identify $H^G_*(pt)$
with $H^*_G$.

Let $\pi^X: X \rightarrow pt$ denote the projection.  If $X$ is
proper, this induces an $H^*_G$-linear map $\pi^X_*: H^G_*(X)
\rightarrow H^G_*(pt) \cong H^*_G$.  In this case, there is a pairing
$(\ , \ ): H^*_G(X) \otimes H^G_*(X) \rightarrow H^*_G$ taking $x
\otimes C$ to $\pi^X_*(x \cap C)$.  We will sometimes write this
pairing as $\int_C x$, and if $C = [Y]_G$, we will abuse notation and
write it as $\int_Y x$.  The pairing has the property that given
$f:X_1 \rightarrow X_2$, we have 
\begin{equation} \label{e.pairing}
(f^*x_2, C_1) = (x_2, f_*C_1).
\end{equation}
(Proof: $(f^*x_2, C_1) = \pi^{X_1}_*(f^* x_2 \cap C_1) =
\pi^{X_2}_*f_*(f^*x_2 \cap C_1) = \pi^{X_2}_*(x_2 \cap f_* C_1) =
(x_2, f_* C_1)$.)  

The map $X \times^G U \rightarrow U/G$ is a fibration with fiber
$X$, and pullback to a fiber yields a map $H^*_G(X) \rightarrow
H^*(X)$.  There is also a Gysin morphism $H_*^G(X) \rightarrow
H_*(X)$.

A variety $X$ is said to be paved by affines if it can be written as a
finite disjoint union $X = \coprod X^0_i$ where $X^0_i$ is isomorphic to
affine space $\A^{d_i}$ for some $d_i$.  As is well known (see
e.g. \cite{Kumar-Nori}) the Borel-Moore homology $H_*(X)$ is the free
$\Z$-module generated by the fundamental classes $[X_i]$ (where $X_i$
is the closure of $X^0_i)$; the odd-dimensional Borel-Moore homology
vanishes.

Part (b) of the next proposition and the remark following are from
\cite{Arabia} (Prop. 2.5.1 and 2.4.1), with a somewhat different
proof.

\begin{proposition} \label{proposition.1}
Suppose the $G$-variety $X$ has a pairing by $G$-invariant affines
$X^0_i$.  Then
\begin{enumerate}
\item[(a)] $H^G_*(X)$ is a free $H^*_G$-module with basis $\{[X_i]_G\}$.
\item[(b)]  Suppose in addition that $X$ is complete and that $H^*_G$
is torsion-free.  Then there exist
classes $x_i$ (of degree $\dim X_i$) in $H^*_G(X)$ which form a basis
for $H^*_G(X)$ as $H^*_G$-module, such that the bases $\{[X_i]_G\}$ and
$\{x_i\}$ are dual in the sense that $\int_{X_i} x_j =
\delta_{ij}$.
\end{enumerate}
\end{proposition}

\noindent{\bf Proof:}  (a)  Let $X^0_k$ be open in $X$, and $Y =
X - X^0_k$; then there is a long exact sequence of
$H^*_G$-modules
$$
\rightarrow H^G_{i+1}(X^0_k) \rightarrow H^G_i(Y) \rightarrow
H^G_i(X)
\rightarrow H^G_i(X^0_k) \rightarrow \cdots
$$
Since $X^0_k$ is isomorphic to affine space, $H^G_*(X^0_k)$ is a free
$H^*_G$-module of rank 1, generated by $[X^0_k]_G$.  Hence all the odd
equivariant homology of $X^0_k$ vanishes, by induction the same holds
for $Y$, and then by the long exact sequence it holds for $X$.  Thus
we have a short exact sequence of $H^*_G$-modules:
$$
0 \rightarrow H^G_*(Y) \rightarrow H^G_*(X) \rightarrow
H^G_*(X^0_k) \rightarrow 0.
$$
This is split by the $H^*_G$-linear map $H^G_* (X^0_k) \rightarrow
H^G_*(X)$ taking $[X^0_k]_G$ to $[X_k]_G$.  Induction implies (a).

(b) Because the odd ordinary cohomology of $X$ vanishes, the pullback
to a fiber $H^*_G(X) \rightarrow H^*(X)$ is surjective (this is
because the spectral sequence of the fibration $X_G \rightarrow B_G$
degenerates at $E_2$).  If $\{y_i\}$ are any classes of pure degree in
$H^*_G(X)$ which pull back to a basis of $H^*(X)$ (we may assume $\deg
y_i = \dim X_i$), then by the Leray-Hirsch theorem \cite{Spanier},
$H^*_G(X)$ is a free $H^*_G$-module with basis $\{y_i\}$.  Claim: The
matrix $A = (a_{ij})$ with entries $a_{ij} \in H^*_G$ defined by
$a_{ij} = \int_{X_i}y_j$ is invertible.  This can be seen by slightly
modifying the arguments of \cite[Theorem 4.1]{Graham}.  For, we may
assume that the $X_i$ are numbered so that the dimension increases as
$i$ increases.  Now, $\int_{X_i} y_j = 0$ unless $\deg y_j \ge \dim
X_i$.  This implies that the matrix $(a_{ij})$ is block upper
triangular (here a block of the matrix corresponds to the set of $(i,
j)$ with $\dim X_i = d$, $\dim X_j = e$, for fixed $d$ and $e$).
Moreover, the diagonal blocks are invertible matrices of scalars (as,
for any fixed $d$, the entries in the corresponding diagonal block are
just the values $([X_i], y'_j)$, where $y'_j$ is the pullback to a
fiber of $y_j$; $\left\{[X_i]\right\}$ is a basis for $H_{2d}(X)$ and
$\{y'_j\}$ a basis for $H^{2d}(X)$).  Hence the matrix $A$ is
invertible, as claimed.

Let $B = A^{-1} = (b_{ij})$ and define $x_j = \sum_i b_{ij} y_i$.  Then
$\{x_j\}$ is a basis of $H^*_G(X)$ dual to $\{[X_i]_G\}$.  Indeed,
\begin{equation} \label{e.determined}
\int_{X_i} x_j = \sum_k \int_{X_i} b_{kj}y_k = \sum_k b_{kj}\int_{X_i}y_k
= \sum_k b_{kj} a_{ik}
= \delta_{ij}.
\end{equation}
Note that the dual basis is uniquely determined by (\ref{e.determined}),
as can be seen by expressing one dual basis in terms of another.
Because the $[X_i]_G$ have pure degree $\dim X_i$, the elements $x_j$
of the dual basis must have degree $\dim X_j$.  For,
if $Y$ is an irreducible closed subvariety of $X$,
and $y \in H^k_G(X)$, then $\int_{Y}y$ has degree $k - \dim Y$.  Hence
if we replace each $x_j$ by its component in degree $\dim X_j$, we
still have a dual basis.  As the dual basis is unique, each
$x_j$ must have degree $\dim X_j$.
\hskip.15in $\square$

\medskip

\textbf{Remarks.} (1) The conditions $\int_{X_i} x_j = \delta_{ij}$ imply that
under the map $H^*_G(X) \rightarrow H^*(X)$, the images $x'_i$ of
$x_i$ form a basis of $H^*(X)$ dual to the basis $\{[X_i]\}$ of
$H_*(X)$.

(2) This result and the next are also valid with coefficients in
a field; then $H^*_G$ is automatically torsion-free.
\medskip

For a variety $X$ paved by $G$-invariant affines as above, we have
the following description of the product
on $H^*_G(X)$ in terms of the diagonal morphism.  The
non-equivariant version of this result was used by \cite{Kumar-Nori}.
The equivariant version was mentioned in \cite{Peterson} for the flag
variety; the general proof is the same.  Note that the diagonal
morphism $\delta: X \rightarrow X \times X$ is $G$-equivariant ($G$
acting diagonally on $X \times X$).

\begin{proposition}  \label{proposition.2}
Let $X$ be a $G$-variety with a paving by $G$-invariant affines
$X^0_i$; assume $H^*_G$ is torsion-free.  Let $X_i$ and $x_i$ be as in
the previous proposition.  We can write $\delta_*[X_k]_G =
\sum_{i,j}a_{ij}[X_i \times X_j]_G$, where $a_{ij} \in H^*_G$.  The
product in $H^*_G(X)$ is given by
$$
x_i x_j = \sum_k a^k_{ij} x_k.
$$
\end{proposition}

\noindent{\bf Proof:}  We can write $\delta_*[X_k]_G$ in the form claimed
because the classes $[X_i \times X_j]_G$ form a basis for $H^G_*(X
\times X)$ as $H^*_G$-module.

Let $q_i: X \times X \rightarrow X$ denote the $i$-th projection.  As
in the non-equivariant case, the product on $H^*_G(X)$ is given by
$$
c_1\cdot c_2 = \delta^*(q^*_1 c_1 \cdot q^*_2 c_2)
$$
for $c_1, c_2 \in H^*_G(X \times X)$.  (This can be seen by considering the
composition
$$
X_G \stackrel{\delta_G}{\rightarrow}(X \times X)_G \cong X_G \times_{BG}X_G
\stackrel{i}{\hookrightarrow} X_G
\times X_G
$$
and noting that the product on $H^*(X_G)$ is given by
$$
\zeta_1 \cdot \zeta_2 = (i \circ \delta_G)^*(pr^*_1 \zeta_1 \cdot pr^*_2 \zeta_2)
$$
where $pr_i: X_G \times X_G \rightarrow X_G$ is the projection and
$\zeta_i \in H^*(X_G)$.  Choosing $\zeta_i$ to represent $c_i \in
H^*_G(X)$, the assertion follows easily.)

The preceding proposition shows that if $X$ is paved by invariant
affines, then $H^*_G(X)$ and $H^G_*(X)$ are free $H^*_G$-modules, with
a perfect pairing
$$
(\ , \ ): H^*_G(X) \otimes_{H^*_G} H^G_*(X) \rightarrow H^*_G .
$$
Using this, we can identify
$$
H^*_G(X) = \mbox{Hom}_{H^*_G}(H^G_*(X), H^*_G) .
$$
Therefore, to show that $x_ix_j = \sum_k a^k_{ij} x_k$, it is enough to
show that for all
$\nu \in H^G_*(X)$, we have
$$
(x_i x_j, \nu) = (\sum_k a^k_{ij} x_k, \nu ) = \sum_k
a^k_{ij}(x_k, \nu).
$$
In fact, it is enough to check this when $\nu$ is one of the basis elements
$[X_k]_G$, i.e., it
is enough to show
$$
(x_ix_j, [X_k]_G)  = a^k_{ij}.
$$

Now
\begin{eqnarray*}
(x_ix_j, [X_k]_G) &=& (\delta^*(q^*_1 x_i \cdot q^*_2 x_j), [X_k]_G)\\
&=& (q^*_1 x_i \cdot q^*_2 x_j, \delta_*[X_k]_G)\\
&=& \sum_{m,n}a^{mn}_k (q^*_1 x_i \cdot q^*_2x_j, [X_m \times X_n]_G).
\end{eqnarray*}
By definition of the pairing, $(q^*_1 x_i \cdot q^*_2 x_j,
[X_m \times X_n]_G) = \pi^{X \times X}_*(q^*_1 x_i \cdot q^*_2 x_j \cap [X_m
\times X_n]_G)$.  This is computed 
using the fibrations $X_G \rightarrow B_G$ and $(X
\times X)_G = X_G \times_{B_G} X_G \stackrel{\pi_G}{\rightarrow}B_G$.  
By the next lemma, the result is equal to
$$
\pi^X_*(x_i \cap [X_m]_G) \cdot \pi^X_*(x_j \cap [X_n]_G)
$$
which is 1 if $i = m$ and $j =n$, and 0 otherwise.  We conclude
$(x_ix_j, [X_k]_G) = a^{ij}_k$, as desired.  \hskip.15in $\square$

\begin{lemma}  \label{lemma.1}
Let $\rho_i: X_i \rightarrow Y$ $(i = 1, 2)$ be fibrations with
$\rho_i$ proper, $\pi: X_1 \times_Y X_2 \rightarrow Y$, $q_i: X_1
\times_Y X_2 \rightarrow X_i$ the projections.  Let $Z_i \subset X_i$
be closed subvarieties and $\alpha_i \in H^*(X_i)$.  Assume $Y$ is
smooth, and identify $H_*(Y)$ with $H^*(Y)$.  Then
$$
\pi_*(q^*_1 \alpha_1 \cdot q^*_2 \alpha_2 \cap [Z_1 \times_YZ_2]) =
\rho_{1*}(\alpha_1 \cap
[Z_1]) \cdot \rho_{2*}(\alpha_2 \cap [Z_2])
$$
where on the right hand side the product is taken in $H^*(Y)$.

\end{lemma}

\noindent{\bf Proof:}  We have a Cartesian diagram
$$
\begin{array}{ccc}
X_1 \times_Y X_2 & \stackrel{\Delta}{\rightarrow} &  X_1 \times X_2 \cr
\downarrow \pi & & \downarrow \Pi \\
Y & \stackrel{\delta}{\rightarrow} & Y \times Y
\end{array}
$$
Because $Y$ is smooth, $\delta$ (and hence $\Delta$) are regular
embeddings, so there are Gysin maps $\delta^*$ and $\Delta^*$ on
homology.  Claim: In $H_*(X_1 \times_Y X_2)$,
$$
q^*_1 \alpha_1 \cdot q^*_2 \alpha_2 \cap [Z_1 \times_Y Z_2] =
\Delta^*((\alpha_1 \cap [Z_1])
\times (\alpha_2 \cap [Z_2])).
$$
To prove this, first note that (with $pr_i: X_1 \times X_2 \rightarrow
X_i$ denoting the projection) $q^*_1 \alpha_1 \cdot q^*_2 \alpha_2=
\Delta^*(pr^*_1 \alpha_1 \cdot pr^*_2 \alpha_2) = \Delta^*(\alpha_1
\times \alpha_2)$ (cf. \cite[p. 351]{Munkres}).  Next, $[Z_1 \times_Y
Z_2] = \Delta^*[Z_1 \times Z_2]$, since $Z_1 \times Z_2$ and $\Delta(X
\times_Y X)$ are subvarieties of $X_1 \times X_2$ whose intersection
at smooth points is transverse.  Hence (noting that $[Z_1 \times Z_2]
= [Z_1]\times [Z_2]$ by [F, p. 377])
$$
\begin{array}{ccc}
q^*_1\alpha_1 \cdot q^*_2\alpha_2 \cap [Z_1 \times_Y Z_2] &=&
\Delta^*(\alpha_1 \times
\alpha_2) \cap \Delta^*[Z_1 \times Z_2]\\
&=& \Delta^*((\alpha_1 \times \alpha_2) \cap ([Z_1]\times [Z_2]))\\
&=& \Delta^*((\alpha_1 \cap [Z_1]) \times (\alpha_2 \cap [Z_2]))
\end{array}
$$
proving the claim.

To complete the proof of the lemma, we compute:
$$
\begin{array}{ccc}
\pi_*(q^*_1 \alpha_1 \cdot q^*_2 \alpha_2 \cap [Z_1 \times_Y Z_2]) &=&
\pi_*\Delta^*((\alpha_1 \cap [Z_1]) \times (\alpha_2 \cap [Z_2]))\cr
&=& \delta^*\Pi_*((\alpha_1 \cap [Z_1]) \times (\alpha_2 \cap [Z_2]))\cr
&=& \delta^*(\rho_{1*}(\alpha_1 \cap [Z_1]) \times \rho_{2*}(\alpha_2 \cap
[Z_2]))\cr
&=& \rho_{1*}(\alpha_1 \cap [Z_1]) \cdot \rho_{2*}(\alpha_2 \cap [Z_2]).
\end{array}
$$
This proves the lemma. \hskip.15in $\square$

\section{The positivity theorem}

In this section we prove the positivity result about multiplication in
equivariant cohomology (Theorem \ref{theorem.1}).  As in the
non-equivariant case considered by Kumar and Nori, it is deduced from
a result about invariant cycles (Theorem \ref{theorem.2}).  In the
non-equivariant setting, Hirschowitz \cite{Hirschowitz} proved that
for a projective scheme with an action of a connected solvable group
$B$, any effective cycle is rationally equivalent to a $B$-invariant
effective cycle.  Kumar and Nori gave a different proof of this result
(without assuming projectivity) in the special case of unipotent groups,
and the proof of Theorem \ref{theorem.2} is adapted from their proof.

In this section, $T$ will denote an algebraic torus (i.e. product of
multiplicative groups $\G_m$) with Lie algebra $\tf = \mbox{Lie }T$,
and $\hat{T} \subset \tf^*$ the group of characters of $T$.  The
equivariant cohomology group $H^*_T$ can be identified with the
polynomial ring $S(\hat{T})$, the symmetric algebra on the
free abelian group $\hat{T}$.

\begin{theorem} \label{theorem.1} 
Let $B$ be a connected solvable group with unipotent radical $N$ and
Levi decomposition $B = TN$.  Let $\alpha_1, \dots, \alpha_d \in
\hat{T}$ denote the weights of $T$ on $\nf = \mbox{Lie }N$.  Let $X$
be a complete $B$-variety on which $N$ acts with finitely many orbits
$X^0_1, \dots, X^0_n$.  These are a paving of $X$ by $B$-stable
affines; let $X_1 \dots, X_n$ denote the closures, so $\{[X_1]_T,
\dots, [X_n]_T\}$ are a basis for $H^T_*(X)$.  Let $\{x_1, \dots,
x_n\}$ denote the dual basis of $H^*_T(X)$.  Write
$$
x_ix_j = \sum_k a^k_{ij} x_k
$$
with $a^k_{ij} \in H^*_T = S(\hat{T})$.  Then each $a^k_{ij}$ can be
written as a sum of monomials $\a^{i_1}_1 \cdots \a^{i_d}_d$, with
nonnegative integer coefficients.
\end{theorem}

Note that the constant term in each $a^k_{ij}$ (i.e., the coefficient
of $\alpha^0_1 \cdots \alpha^0_d)$ is nonnegative by the above
theorem.  This is the coefficient that occurs in the multiplication
in the ordinary cohomology $H^*(X)$.  The reason is that our
hypotheses imply $H^*(X) = H^*_T(X)/H^{>0}_T \cdot H^*_T(X)$
(see \cite{GKM}).

The next result is the key ingredient in the proof of Theorem
\ref{theorem.1}.  In this theorem, $N$ is not assumed to act with 
finitely many orbits.  The result also holds with equivariant Chow
groups in place of equivariant Borel-Moore homology.

\begin{theorem}  \label{theorem.2}
Let $B$ be a connected solvable group with unipotent radical $N$, and
let $T \subset B$ be a maximal torus, so $B = TN$.  Let $\alpha_1,
\dots, \alpha_d \in \hat{T}$ denote the weights of $T$ acting on 
$\nf = \mbox{Lie }N$.  Let $X$ be a scheme with a $B$-action and $Y$ a
$T$-stable subvariety of $X$.  Then there exist $B$-stable
subvarieties $D_1, \dots, D_r$ of $X$ such that in $H^T_*(X)$,
$$
[Y]_T = \sum f_i [D_i]_T
$$
where each $f_i \in H^*_T$ can be written as a linear combination of
monomials in $\alpha_1, \dots, \alpha_d$ with nonnegative integer
coefficients.
\end{theorem}

The following lemma was pointed out to me by Michel Brion.

\begin{lemma} \label{l.orbits}
Suppose the connected solvable group $B = TN$ acts on $X$ and that
$N$ has finitely many orbits on $X$.  Then each $N$-orbit is $B$-stable
(in fact, the $B$-orbit of a $T$-fixed point).
\end{lemma}

\medskip

\noindent{\bf Proof:} $B$ has finitely many orbits on $X$ (as the subgroup $N$
does); as each $N$-orbit is $N$-stable, it is a finite union of
$N$-orbits.  Let $B \cdot x' \simeq B/B'$ be an orbit, where $B'$ is
the stabilizer of $x'$.  As each $N$-orbit is isomorphic to affine
space (see e.g. \cite{Kumar-Nori}), the odd cohomology of $B \cdot x'$
vanishes, so $B'$ must contain a maximal torus of $B$.  As all maximal
tori of $B$ are conjugate \cite[Corollary 11.3]{Borel}, there is some
$b \in B$ such that $B' = b B_1 b^{-1}$, where $B_1 \supset T$.  Then
$B \cdot x' = B \cdot x$ where $x = b^{-1} x'$; moreover $B_1$ is the
stabilizer of $x$.  Hence $B \cdot x$ is the $N$-orbit of the
$T$-fixed point $x$. \hskip.15in $\square$

\medskip

\noindent{\bf Proof of Theorem \ref{theorem.1}:} The group $\tilde{B}
= T \cdot (N \times N)$ (semi-direct product) acts on $X \times X$ by
$t \cdot (n_1, n_2) (p_1, p_2) = (tn_1 p_1, tn_2 p_2)$.  The unipotent
radical $N \times N$ has finitely many orbits $X^0_i \times X^0_j$ on
$X \times X$, with closures $X_i \times X_j$, so $H^T_*(X \times X)$
is a free $H^*_T$-module with basis $[X_i \times X_j]_T$.  By
Proposition \ref{proposition.2}, if $x_i x_j = \sum_k a^k_{ij} x_k$
then $\delta_*[X_k]_T = [\delta(X_k)]_T= \sum_{ij} a^k_{ij}[X_i \times
X_j]_T$. The coefficients $a^k_{ij}$ are uniquely determined by the
expansion of $\delta_*[X_k]_T$ because the classes $[X_i \times X_j]$
are linearly independent over $H^*_T$.  By Theorem \ref{theorem.2},
these coefficients can be written as monomials in $\alpha_1, \dots,
\alpha_d$ with nonnegative integer coefficients, where $\alpha_1,
\dots, \alpha_d$ are the weights of $T$ on $\mbox{Lie }(N \times N)$
(which are the same as the weights of $T$ on $\nf$).  \hskip.15in
$\square$

\bigskip

\noindent{\bf Proof of Theorem \ref{theorem.2}:} First consider the
case where $\dim N = 1$; then $B/T \stackrel{\sim}{\rightarrow} N
\stackrel{\varphi}{\rightarrow} \G_a$, where $\G_a \cong \A^1$ is the
additive group.  Write $\alpha = \alpha_1$.  We have $B = NT$, and the
map $B/T \stackrel{\sim}{\rightarrow} N$ sends $n T \rightarrow n$.
Now, $B$ acts on $B/T$ by left multiplication.  Via the isomorphism
of $B/T$ with $N$, we obtain an action of $B$ on $N$; the subgroup $T
\subset B$ acts on $N$ by conjugation, and the subgroup $N$ acts by
left multiplication.  The action of $T$ by conjugation on $N$
corresponds under $\varphi$ to an action of $T$ on $\A^1$ with weight
$\alpha$.  Embed $B/T \hookrightarrow \P^1$ by $n T \mapsto
[\varphi(n): 1]$.  The action of $B$ on $B/T$ extends to an action on
$\P^1$: the element $tn \in B$ acts by the matrix
$$
\left( \begin{array}{cc}
      \alpha(t) & \varphi(n) \\
       0 & 1
\end{array} \right).
$$
The point $\infty = [1:0]$ is fixed by $B$, while the point $0
= [0:1]$ is fixed by $T$.

Now, $B$ acts on $B \times^T X$ by left multiplication: $b \cdot (b',
x) =(bb', x)$.  Under the isomorphism $\theta: B \times^T X
\rightarrow B/T \times X$ taking $(b, x)$ to $(bT, bx)$, the
$B$-action corresponds to the product action on $B/T \times X$.  This
extends to a $B$-action on $\P^1 \times X$.  The projections $\pi: \P^1
\times X \rightarrow \P^1$ and $\rho: \P^1 \times X \rightarrow X$ are
$B$-equivariant.

If $Y \subset X$ is a $T$-invariant subvariety then $B \times^T Y$ is
a $B$-invariant subvariety of $B \times^T X$.  Let $Z$ be the Zariski
closure of $\theta(B \times^T Y)$ in $\P^1 \times X$; $\theta(B
\times^T Y)$ and $Z$ are $B$-invariant subvarieties of $\P^1 \times X$.
Let $\pi_Z$ denote the restriction of $\pi$ to $Z$.

Let $[w_0: w_1]$ be projective coordinates on $\P^1$, and $w$ the
rational function $\frac{w_0}{w_1}$.  Let $g = \pi^*_Z w$; then $w$
(and hence $g$) are rational functions which are $T$-eigenvectors of
weight $-\alpha$.  By \cite[Theorem 2.1]{Brion}\footnote{Brion is
using the convention that if $X$ is a $T$-space, then $T$ acts on
functions on $X$ by $(t\cdot f)(x) = f(tx)$, while we are using the
convention that $T$ acts on functions by $(t\cdot f)(x) = f(t^{-1}
x)$.  Under Brion's convention, our function $g$ would be an
eigenvector of weight $\alpha$.} we have in $H^T_*(\P^1 \times X)$ the
relation $[\mbox{div}_Z g]_T = \alpha[Z]_T$.  Therefore, in $H^T_*X$ we
have the relation
\begin{equation} \label{e.11}
\rho_{*}[\mbox{div}_Z g]_T = \alpha \rho_{*}[Z]_T
\end{equation}

Now, $\pi^{-1}_Z(0) = \{0\} \times Y$ (cf. \cite{Kumar-Nori}).  Also,
$\pi^{-1}_Z (\infty) = \{\infty\}\times D$ where $D$ is a subscheme of
$X$.  Therefore (\ref{e.11}) yields
$$
[Y]_T = [D]_T + \alpha \rho_{*}[Z]_T
$$
As $\pi_Z$ is $B$-equivariant, and $\infty \in \P^1$ is $B$-fixed, it
follows that $\{\infty\} \times D$, and hence $D$, are $B$-invariant.
Each irreducible component $D_i$ ($i = 1,\dots, r)$ of $D$ is
therefore $B$-invariant (as $B$ is connected) and if $m_i$ is the
multiplicity of $D_i$ in $D$ then $[D]_T = \sum^r_{i=1}m_i[D_i]_T$.
Likewise, $\rho$ is $B$-equivariant and $Z$ is $B$-invariant.  If
$Z_i$ is a component of $Z$ then the map $\rho|_{Z_i}$ of $Z_i$ onto
its image in $X$ is finite if and only if the map $\rho_T|_{Z_{iT}}$
of $Z_{iT}$ onto its image in $X_T$ is finite, and the degrees of the
maps are the same.  If we list the components of $\rho(Z)$ which are
finite images of components of $Z$ as $D_{r+1}, \dots, D_s$, it
follows that each of these components is $B$-invariant and that
$\rho_{*}[Z]_T = \sum^s_{i=r+1} m_i[D_i]_T$ where $m_i$ are positive
integers.  We conclude that
\begin{equation} \label{e.22}
[Y]_T = \sum^r_{i=1}m_i[D_i]_T + \sum^s_{i=r+1}m_j\alpha[D_i]_T
\end{equation}
where the $D_i$ are $B$-invariant.  This proves the result if $\dim N= 1$.

To prove the result in general, we can find a subgroup $N' \subset N$
such that $N'$ is normal in $B$ and $\dim N/N' = 1$.  Let $\alpha$ be
the weight of $T$ on $\mbox{Lie }(N/N')$.  Define $B' = N'T \subset B
= NT$.  By induction, we may assume the result is true for $B'$.  It
is enough to show that given a $B'$-invariant subvariety $Y \subset
X$, we can write $[Y]_T$ as in (\ref{e.22}), with $B$-invariant $D_i$.
For this we modify the above proof, as follows.  Replace $B/T$,
$B\times^T X$, and $B\times^T Y$ by $B/B'$, $B\times^{B'}X$, and
$B\times^{B'} Y$; the map $\theta$ now takes $B\times^{B'}X$ to
$B/B'\times X$.  Again $\varphi: B/B' \stackrel{\cong}{\rightarrow}
\G_a = \A^1$ and $T$ acts by weight $\alpha$ on $\A^1$.  We can embed
$B/B' \hookrightarrow \P^1$ as before; the point $\infty = [1:0]$ is
fixed by $B$, and $[0:1]$ is fixed by $B'$.  With these modifications,
(\ref{e.22}) is proved as above.  This proves the theorem.
\hskip.15in $\square$

\section{Schubert varieties}

\subsection{Peterson's conjecture}

Let $G$ be a complex semisimple group and $B \supset T$ a Borel
subgroup and maximal torus, respectively.  Let $N$ be the unipotent
radical of $B$; let $B^- = TN^-$ be the opposite Borel.  Choose a
system of positive roots so that the roots in $\nf$ are positive.  Let
$W = N(T)/T$ denote the Weyl group; we abuse notation and write $w$
for an element of $W$ and also for a representative in $N(T)$.  Let $X
= G/B$ the flag variety.  The $T$-fixed points are $\{wB\}_{w\in W}$;
let $X^0_w = N \cdot wB \subset X$ and $Y^0_w = N^- \cdot wB$.  Then
$X = \coprod_w X^0_w$ (resp $X = \coprod_w Y^0_w$) is a decomposition
of $X$ as a disjoint union of finitely many $N$ (resp.  $N^-$)-orbits.
Let $X_w$ and $Y_w$ denote the closures of $X^0_w$ and $Y^0_w$, and
$\{x_w\}$ and $\{y_w\}$ the bases of $H^*_T X$ dual (in the sense of
Proposition \ref{proposition.1}) to $\{[X_w]_T\}$ and $\{[Y_w]_T\}$.

Let $\alpha_1, \dots, \alpha_{\ell}$ denote the simple roots.  Any
weight of $T$ on $\nf$ (resp. $\nf^-$) is a nonnegative
(resp. nonpositive) linear combination of the simple roots.
Therefore, the next corollary is an immediate consequence of Theorem
\ref{theorem.1}.

\begin{corollary}  \label{corollary.1}
With notation as above, write $x_u x_v = \sum_w a^w_{uv}
x_w$ and $y_u y_v = \sum_v b^w_{uv} y_w$, with $a^w_{uv}$ and
$b^w_{uv}$ in $H^*_T$.  Then $a^w_{uv}$ (resp.  $b^w_{uv}$) is a
linear combination of monomials in the $\alpha_i$, with nonnegative
(resp.  nonpositive) coefficients.  \hskip.15in $\square$
\end{corollary}

\remark Theorem \ref{theorem.1} can be applied to the varieties $X_w$
and $Y_w$, which are in general singular, to yield an analogue of
Corollary \ref{corollary.1} for $H^*_T(X_w)$ and $H^*_T(Y_w)$.  The
analogous result also holds for partial flag varieties.

\medskip

Because $X$ is smooth, the map $H^*_T(X) \stackrel{\cap
[X]_T}{\rightarrow} H^T_*(X)$ is an isomorphism.  The next lemma
is known (cf. \cite{Peterson}) but for lack of reference we give a proof.

\begin{lemma}  The map $H^*_T(X) \stackrel{\cap [X]_T}{\rightarrow}
H^T_*(X)$ takes $y_w$ to $[X_w]_T$.
\end{lemma}

\noindent{\bf Proof:} We can identify $H^*_T(X)$ with
$\mbox{Hom}_{H^*_T}(H^T_*(X), H^*_T)$ (see the proof of Proposition
2). Hence, any $\gamma \in H^*_T(X)$ is uniquely determined by the
values $\pi^T_*(\gamma \cap h')$ as $h'$ ranges over the basis
$\{[Y_{w'}]_T\}$ of $H^T_*(X)$.

Now, if $\gamma \in H^*_T(X)$ satisfies $\gamma \cap [X]_T = h$, then
$\gamma \cap h' = h \cdot h'$.  Indeed, the intersection product on
$H^T_*(X)$ satisfies: if $\gamma' \cap [X]_T = h'$, then $\gamma \cdot
\gamma' \cap [X]_T = h \cdot h'$; but $\gamma \cdot \gamma' \cap
[X]_T = \gamma \cap (\gamma' \cap [X]_T) = \gamma \cap h'$.

Combining these facts, we see that to show $y_w \cap [X]_T = [X_w]_T$,
it suffices to show
$$
\pi^X_*([X_w]_T[Y_{w'}]_T) = \pi^X_*(y_w \cap [Y_{w'}]_T) = \delta_{ww'}.
$$
Now, for any $w, w'$, the intersection $X_w \cap Y_{w'}$ is
$T$-invariant, and is known to satisfy $\mbox{codim }X_w \cap Y_{w'} =
\mbox{codim }X_w + \mbox{codim }Y_{w'}$.  (Indeed, by
\cite{Kazhdan-Lusztig}, $X_w \cap Y_{w'}^0$ is irreducible and of
dimension $\dim X - \dim X_w - \dim Y_{w'}$, but by
\cite[p. 137]{Fulton}, each component of $X_w \cap Y_{w'}$ has at
least that dimension.  It follows that $X_w \cap Y_{w'}^0$ is dense in
$X_w \cap Y_{w'}$.)  Hence $[X_w]_T[Y_{w'}]_T$ is a multiple of $[X_w
\cap Y_{w'}]_T$.  If $\dim X_w \cap Y_{w'} > 0$, then $\dim(X_w \cap
Y_{w'})_T > \dim BT$, so $\pi^X_*([X_w \cap Y_{w'}]_T) = 0$.  If $\dim
X_w \cap Y_{w'} = 0$, then $w = w'$ and $X_w$ and $Y_w$ intersect with
multiplicity $1$ at the point $wB$ \cite[Prop. 2]{Chevalley}.  Hence
$\pi^X_T:X_T \rightarrow BT$ maps $(X_w \cap Y_w)_T$ isomorphically
onto $BT$, and therefore $\pi^X_*([X_w]_T[Y_w]_T)= \pi^X_*([X_w\cap
Y_w]_T) = 1$.  This proves the lemma. \hskip.15in $\square$
\medskip

The intersection product on $H^T_*(X)$ is induced by the product on
$H^*_T(X)$, via the isomorphism $\cap [X]_T$. The above lemma and
Corollary \ref{corollary.1} therefore imply:

\begin{corollary} \label{corollary.2}
The intersection product on $H^T_*(X)$ is given by $[Y_u]_T [Y_v]_T =
\sum_w a^w_{uv}[Y_w]_T$ (resp. $[X_u]_T [X_v]_T =
\sum_w b^w_{uv}[X_w]_T$), where each $a^w_{uv}$ (resp. $b^w_{uv}$) in $H^*_T$ 
is a sum of monomials in the $\alpha_1, \dots, \alpha_{\ell}$, with
nonnegative (resp. nonpositive) coefficients.\hskip.15in $\square$
\end{corollary}

Corollaries \ref{corollary.1} and \ref{corollary.2} were conjectured by
Dale Peterson.

\medskip

\textbf{Example.}  As a concrete example, we work out the case of the flag
variety of $SL_2$.  Here, we take $B$ (resp. $B^-$, $T$) to be the
upper triangular (resp. lower triangular, diagonal) matrices; we
identify $X$ with $\P^1$, acting as usual.  Then $W = \{1, s\}$ and
the Schubert varieties are $X_1 = [1: 0]$, $X_s = X$, $Y_1 = X$, $Y_s
= [0:1]$.  The character group of $T$ is $\hat{T} = \Z \cdot x \cong
\Z$, and the positive root is $\alpha = 2x$.  The ring $H^*_T = \C[x]$.  
The action of $T$ on $\P^1$ is with weights $\pm 1$, so $H^*_TX =
\C[x, h]/(h + x)(h - x)$.  We will identify $H^*_T X$ with $H^T_* X$
via $\cap [X]_T$.  Under this isomorphism, $[X_s]_T = [Y_1]_T = 1$.
If $[z_0: z_1]$ are projective coordinates on $\P^1$, then $z_0$ may
be viewed as a section of $\O(1)$ which is a $T$-eigenvector of weight
$-1$.  Then $z_0 \otimes 1$ is a $T$-invariant section of $\O(1)
\otimes \C_1$ (here $\C_1$ is the trivial line bundle with $T$ with
weight 1).  The zero-scheme of $z_0 \otimes 1$ is $[0:1]$, so we
conclude $[Y_s]_T = [0:1]_T = c^T_1(\O(1)\otimes \C_1) = h + x$.
Similarly, $[X_1]_T = h - x$.  The only interesting multiplication
among the classes $[X_w]_T$ is
$$
[X_1]_T \cdot [X_1]_T = (h - x)^2 = h^2 - 2hx + x^2 = 2x^2 - 2hx =
-2x(h - x) = -\alpha[X_1]_T.
$$
Similarly, the only interesting multiplication among the classes
$[Y_w]_T$ is
$$
[Y_s]_T[Y_s]_T = \alpha[Y_s]_T.
$$
These agree with Corollaries \ref{corollary.1} and \ref{corollary.2}.

\subsection{Billey's conjecture}
Kostant and Kumar \cite{Kostant-Kumar} defined functions (for each $w
\in W$) $\xi^w: W \rightarrow S(\hat{T}) \subset S(\tf^*)$, and showed
that for any $u,v \in W$, one can write
$$
\xi^u \xi^v = \sum_w p^{uv}_w \xi^w
$$
for unique $p^{uv}_w \in S(\tf^*)$.  Billey \cite{Billey} observed in
examples that that if $\nu \in \tf$ satisfies $\alpha(\nu) > 0$ for
all positive roots $\alpha$, then $p^{uv}_w(\nu) \geq 0$, and asked
if a geometric proof was possible.

Arabia \cite{Arabia} proved the following relation of the functions
$\xi^w$ to the $T$-equivariant equivariant cohomology of the flag
variety.  We use the notation of the preceding subsection: thus, $i_w:
wB \rightarrow G/B = X$ denotes the inclusion, and $i_w^*: H^*_T(X)
\rightarrow H^*_T(wB) = H^*_T$ the pullback.  As usual, we identify
$H^*_T(X)$ with $H_*^T(X)$.  

\begin{theorem} \label{p.kk}
\begin{enumerate}
\item[(1)] $i_u^*x_w =\xi^{w^{-1}}(u^{-1})$.
\item[(2)] $p_{w^{-1}}^{u^{-1},v^{-1}} = a_{uv}^w$.
\end{enumerate}
\end{theorem}
\medskip
This is proved (in the general Ka\v{c}-Moody case) in \cite[Theorem
4.2.1]{Arabia}.  We have stated this theorem using the conventions of
\cite{Kostant-Kumar} for the functions $\xi^w$; below we explain the
relationship between the conventions of \cite{Arabia} and
\cite{Kostant-Kumar}.  Note that (2) follows immediately from (1),
since (as noted by Arabia) the pullback $\oplus i_w^*: H^*_T(X)
\rightarrow \oplus H_*^T$ is injective.

As a consequence, we obtain Billey's conjecture:

\begin{corollary} \label{c.billey}
If $\nu \in \tf$ satisfies $\alpha(\nu) > 0$ for all positive roots
$\alpha$, then $p^{uv}_w(\nu) \geq 0$.
\end{corollary}

\noindent{\bf Proof:} This follows immediately from the preceding corollary
and Corollary \ref{corollary.2}.  \hskip.15in $\square$

\medskip

We now discuss the conventions of \cite{Arabia} and
\cite{Kostant-Kumar}.  Let $\C[W]$ denote the group algebra over $\C$
of $W$; let $Q$ be the quotient field of
$S(\tf^*)$.  Kostant and Kumar set $Q_W = \C[W] \otimes Q$; Arabia
defines $Q$ and $Q_W$ with rational rather than complex coefficients,
but we will ignore this difference.  Both \cite{Kostant-Kumar} and
\cite{Arabia} define elements $\xi^w \in \mbox{Hom}_{Q}(Q_W,Q)$, but
with different conventions: if we use $\xi^w$ for the elements defined
in \cite{Kostant-Kumar} and $\xi^w_A$ for the elements defined in
\cite{Arabia}, then $\xi^w = \xi^{w^{-1}}_A$.

Let $F(W,Q)$ denote the set of functions from $W$ to $Q$.  Both
\cite{Kostant-Kumar} and \cite{Arabia} use identifications $F(W,Q)
\stackrel{\simeq}{\rightarrow} \mbox{Hom}_{Q}(Q_W,Q)$; we will denote
their respective identifications by
$$
\begin{array}{ccc}
f \mapsto f_K, & \mbox{where } f_K(\delta_u \otimes 1) = f(u) &
\cite[(4.17)]{Kostant-Kumar} \\
f \mapsto f_A, & \mbox{where } f_A(\delta_u \otimes 1) = f(u^{-1}) &
\cite[\mbox{Section 4.1}]{Arabia}.
\end{array}
$$
If we define $f^w$ and $g^w$ in $F(W,Q)$ by $f^w_K = \xi^w$, 
$g^w_A = \xi^w_A$, then $f^w(u) = g^{w^{-1}}(u^{-1})$.

Arabia uses the injection
$$
\oplus i^*_u: H^*_T(X) \hookrightarrow \oplus H^*_T \simeq F(W,
S(\tf^*)) \subset F(W,Q)
$$
to identify $H^*_T(X)$ with a subset of $F(W,Q)$.  In his paper, he
proves that under this identification, $g^w$ corresponds to what we
have denoted by $x_w \in H^*_T(X)$.  In
\cite{Kostant-Kumar} there is no separate notation introduced for the
$f^w$, but rather they are identified with $\xi^w$, i.e., the $\xi^w$
are viewed as elements of $F(W,Q)$.  If we return to their notation,
we see $\xi^{w^{-1}}({u^{-1}}) = i^*_u x_w$, as stated in Theorem
\ref{p.kk}.

Note that if we let $\xi^w_B$ denote the functions used by Billey,
then $\xi^w_B(u) = \xi^{w^{-1}}(u^{-1})$.

\subsection{The Ka\v{c}-Moody case}
The analogues of Corollaries \ref{corollary.1} and \ref{c.billey} are
also valid for flag varieties (complete or partial) of Ka\v{c}-Moody
groups.  The key point is that such a flag variety, although in
general infinite dimensional, can be approximated by finite
dimensional varieties for which the hypotheses of Theorem
\ref{theorem.1} are satisfied.  Indeed, this was exactly the geometric
motivation of Kumar and Nori.  We will briefly sketch how this works
in equivariant cohomology.  The basic facts we need can be found in
\cite{Slodowy}, to which we refer for a more detailed explanation of
the notation.  Let $G$ be a Ka\v{c}-Moody group and $B$ a Borel
subgroup; let $X = G/B$ denote the flag variety.  The group $B$ is a
proalgebraic group (inverse limit of algebraic groups), and it has a
Levi decomposition $B = T N$, where $N$ is a proalgebraic prounipotent
group (denoted by $U$ in \cite{Slodowy} and \cite{Kumar-Nori}) and $T$
is a finite dimensional torus.  The space $X$ has the structure of
ind-variety: it is realized as a union $X = \cup_{k>0}X_k$, where each
$X_k$ is a finite dimensional variety embedded as a closed subvariety
of $X_{k+1}$.  Here $X_k$ is defined as follows.  We have $X = \coprod
X_w^0$, realizing $X$ as a disjoint union of Schubert cells $X_w^0 = B
\cdot w B$.  The union is over all elements of the Weyl group $W$;
each $X_w^0$ is isomorphic to the affine space $\A^{l(w)}$, where
$l(w)$ is the length of $w$.  By definition, $X_k = \coprod_{l(w) \leq
k} X_w^0$; this is a finite dimensional projective variety which is
paved by affines.  Moreover, each $X_k$ is $B$-stable, and there
exists a subgroup $N_k \subset N$, normal in $B$, such that $B_k =
B/N_k$ is a finite dimensional solvable group, and the action of $B$
on $X_k$ factors through the map $B \rightarrow B_k$.  Each $X_k$
therefore satisfies the hypotheses of Theorem \ref{theorem.1}.  As in
the finite case, there is a set of simple roots $\alpha_1, \ldots,
\alpha_l$ in $\tf^*$, and moreover, for any $k$, every weight in
$\mbox{Lie }(N/N_k)$ is a nonnegative linear combination of simple
roots.

Now, for any fixed $i$, the pullback $H^i_T(X) \rightarrow H^i_T(X_k)$ is a
canonical isomorphism for $k$ sufficiently large (as the decomposition
of $X$ into Schubert cells makes $X$ a CW-complex, and $X_k$ contains all
cells in $X$ of dimension $\leq 2k$, and similarly for the mixed spaces
$X_{kT}$ and $X_T$).  There is a basis $\{x_w \}$ of $H^*_T(X)$ dual
to the fundamental classes $[X_w]_T$, in the sense that the pullbacks
to $H^*_T(X_k)$ form a basis dual to the $[X_w]_T \in H^T_*(X_k)$, for
$l(w) \leq k$.  This basis does not depend on $k$, as can be seen using
property (\ref{e.pairing}) of the pairing, applied to the inclusion map
of $X_k$ into $X_{k+1}$.  Theorem \ref{theorem.1} therefore implies the
following corollary, also conjectured by Peterson.

\begin{corollary}  \label{corollary.3}
With notation as above, if $X$ is the flag variety of a Ka\v{c}-Moody
group, with basis $\{x_w \}$ of $H^*_T(X)$, then $x_u x_v = \sum_w a^w_{uv}
x_w$, with $a^w_{uv} \in H^*_T$ a
linear combination of monomials in the $\alpha_i$, with nonnegative
coefficients.  \hskip.15in $\square$
\end{corollary}

\end{document}